\documentclass{svmult-ddm}

\usepackage{mathptmx}       
\usepackage{helvet}         
\usepackage{courier}        
\usepackage{type1cm}        
\usepackage{graphicx}        

\usepackage[bottom]{footmisc}
\usepackage{amssymb}
\usepackage{amsfonts}
\usepackage{amsbsy}
\usepackage{amscd}
\usepackage{amstext}
\usepackage{dsfont}
\usepackage[english]{babel}
\usepackage{color}
\usepackage{graphics}
\usepackage{wrapfig}
\usepackage{psfrag}
\usepackage{color}
\usepackage{url}
\usepackage{verbatim}
\usepackage{algorithm}
\usepackage{algorithmic}

\begin{document}

\title*{Space Decompositions and Solvers for Discontinuous Galerkin
  Methods\thanks{Submitted to the Proceedings of the 21st International Conference on Domain Decomposition Methods}}


\author{Blanca Ayuso de Dios\inst{1}\and Ludmil Zikatanov\inst{2}}

\institute{Centre de Recerca Matem\`atica, Barcelona,
  Spain. \texttt{bayuso@crm.cat}
  \and Department of Mathematics, The Pennsylvania State University,
  University Park, USA
\texttt{ludmil@psu.edu}
}

\maketitle

\abstract{We present a brief overview of the different domain and space
  decomposition techniques that enter in developing and analyzing
  solvers for discontinuous Galerkin methods. Emphasis is given to the
  novel and distinct features that arise when considering DG discretizations over
  conforming methods. Connections and differences with the conforming
  approaches are emphasized.}

\keywords{discontinuous Galerkin, Schwarz methods, Space decompositions}

\section{Introduction}\label{ayusob_plenary_sec:1}

The design and the analysis of efficient preconditioners for
discontinuous Galerkin discretizations has been subject of intensive
research in the last decade with efforts focused mainly on elliptic
problems. 

A standard point of view when studying most of the preconditioning and
iterative solution strategies, in general, is associated with a
particular {\it space decomposition}. From the classical theory of
Lions~\cite{ayusob_plenary_lions1,ayusob_plenary_toselli-widlund,ayusob_plenary_XZ},
we know that, the choice of the space decomposition plays 
significant role in the construction and also in the
convergence properties of the resulting preconditioners. For
nonconfoming methods, domain decomposition and multigrid
preconditioners have been analyzed by establishing connections with
their respective conforming sub-spaces
\cite{ayusob_plenary_brennerNC,ayusob_plenary_Oswald96}.  In the case
of DG methods, the discontinuous nature of the DG finite element
spaces allows to introduce and study not only space splittings
pertinent to the conforming methods but also consider new splittings
which give rise to new techniques and ideas.

In most of the earlier works, relevant space splittings of the DG
finite element space, were introduced via a domain decomposition.
Overlapping additive Schwarz methods have been studied following the
classical Schwarz theory for different DG schemes
\cite{ayusob_plenary_kara0,ayusob_plenary_bren10,ayusob_plenary_wid0}. Contrary
to the conforming case, additive (and multiplicative) Schwarz methods
based on non-overlapping decomposition of the computational domain
have been constructed and proven to be convergent for DG methods. For
such type of preconditioners, novel features, which have no analog in
the conforming case, arise. 
For both overlapping and non-overlapping Schwarz methods, the
splittings are stable in the $L^{2}$-norm by construction and can be
shown to be stable in the natural DG energy norm, with constants
depending on the mesh sizes relative to the coarse and fine subspaces.

More sophisticated substructuring preconditioners have been studied
recently for two dimensional elliptic Poisson problems. In
\cite{ayusob_plenary_sarkis0,ayusob_plenary_sarkis1,ayusob_plenary_sarkis2,
  ayusob_plenary_noi} non-overlapping BDDC, N-N, FETI-DP and BPS
domain decomposition preconditioners are introduced and analyzed for a
Nitsche-type approximation. BDDC preconditioners are studied in
\cite{ayusob_plenary_luca, ayusob_plenary_Joachim} for IP-spectral and IP-hybridized methods. Also here, several
different approaches have been considered and new theoretical tools
have been introduced. And of course, the space splitting in which the
preconditioner rely, comes always from domain decomposition.
Starting directly with a splitting of the DG space, dictated by a
hierarchy of meshes, multigrid methods have been proposed and analyzed
in \cite{ayusob_plenary_GK,ayusob_plenary_brennerMGIP}.  A different
approach was taken in \cite{ayusob_plenary_dobrev} and
\cite{ayusob_plenary_dahmen1,ayusob_plenary_kolja1}, to develop
respectively, two-level and multilevel preconditioners for the
Interior Penalty (IP) DG methods. A common idea behind these works is
to use the fictitious/auxiliary spaces for which one knows how do
develop a preconditioner. Such preconditioning techniques have already
been applied in a wide range of problems in the conforming case. 

The aforementioned auxiliary space preconditioners use error
corrections from the conforming finite element space and they are
certainly related to the a posteriori theory for DG
methods~\cite{ayusob_plenary_kara-apost1}.  In fact, the stable projections given in
\cite{ayusob_plenary_kara-apost1} provide the required tools for constructing
and analyzing the convergence of these preconditioners including the
case of non-conforming meshes.

A novel approach was taken in \cite{ayusob_plenary_az} where a natural
decomposition of the linear DG finite element space was
introduced. The components of the space decomposition are orthogonal
in the inner product provided by the DG bilinear form.  Such a
splitting allows to devise efficient multilevel methods and uniform
preconditioners and analyze these iterative schemes in a clean and
transparent way. This seems to be the only approach available
till now, to prove convergence for the solvers of the {\it
  non-symmetric} Interior Penalty methods. While the methodology has
been applied to the lowest order DG space and conforming meshes, it
is valid in two and three dimensions, and has already been adapted and
extended to a larger family of problems: elliptic
with jump coefficients \cite{ayusob_plenary_jump}; linear elasticity
\cite{ayusob_plenary_AyusoB_GeorgievI_KrausJ_ZikatanovL-2009aa}; and
convection dominated problems corresponding to drift-diffusion models
for transport of species \cite{ayusob_plenary_ariel00}.

We present here a brief overview of some of the domain and space
decomposition techniques that comprise a set of key tools used in
developing and analyzing solvers for DG methods.
In Section \ref{ayusob_plenary_sec:3} we focus on non-overlapping
Schwarz domain decomposition methods.  In Section
\ref{ayusob_plenary_sec:4} and \ref{ayusob_plenary_sec:5} we present
the two main classes of space decomposition methods commenting on
their strengths and weaknesses.

\section{Discontinuous Galerkin Methods}\label{ayusob_plenary_sec:2}
We consider the model problem  for given data  $f\in L^2(\Omega)$:
\begin{equation}\label{ayusob_plenary_mod0}
-\Delta u^{\ast}=f \quad\textrm{in $\Omega$}\qquad \quad
u^{\ast}=0 \quad  \textrm{on $\partial \Omega$}\;,
\end {equation}
Here, $\Omega \subset \mathbb{R}^{d}$, $d=2,3$ is a polygonal
(polyhedral) domain. Let $\mathcal{T}_h$ be a shape-regular family of
partitions of $\Omega$ into $d$-dimensional simplexes $T$ (triangles
if $d=2$ and tetrahedrons if $d=3$) and let $h=\max_{T \in
  \mathcal{T}_h} h_{T}$ with $h_{T}$ denoting the diameter of $T$ for
each $T \in \mathcal{T}_h$.  We denote by $\mathcal{E}^{o}_h$ and
$\mathcal{E}^{\partial}_h$ the sets of all interior faces and boundary
faces (edges in $d=2$), respectively, and we set
$\mathcal{E}_h=\mathcal{E}^{o}_h\cup \mathcal{E}^{\partial}_h$.  Let
$V_h^{DG}$ denote the discontinuous finite element space defined by:
\begin{equation} 
V_{h}^{DG}=\left\{  u\in L^{2}(\Omega) \,\, : \,\, u_{|_{T}}\,\,  \in\,\,  \mathbb{P}^{\ell}(T) \,\, \forall  T \in \mathcal{T}_h\,\right\},\label{ayusob_plenary_defDG} 
\end{equation}
where $\mathbb{P}^{\ell}(T)$ denotes the space of  polynomials of degree at most $\ell$ on each $T$. We also define the conforming finite element space as   
 $V_{h}^{{\rm conf}}=V_{h}^{DG}\cap H^{1}_0(\Omega)$.\\
We define the \emph{average} and \emph{jump} trace operators. Let $T^{+}$ and $T^{-}$ be two neighboring elements, and ${\bf n}^{+}$, ${\bf n}^{-}$ be their outward normal unit vectors, respectively (${\bf n}^{\pm}={\bf n}_{T^{\pm}}$). Let
$\zeta^{\pm}$ and ${\boldsymbol {\tau}}^{\pm}$ be the restriction of $\zeta$ and
${\boldsymbol {\tau}}$ to $T^{\pm}$. We set:
\vskip -0.1cm 
\begin{equation}\label{ayusob_plenary_av-jump}
\begin{array}{ccccc}
2\{\zeta\}&=(\zeta^+ +\zeta^-),\quad
{\lbrack\!\lbrack\, \zeta \,\rbrack\!\rbrack}
=\zeta^+{\bf n}^++\zeta^-{\bf n}^- \quad &\mbox{on } E\in
\mathcal{E}^{o}_h,\\
2\{{\boldsymbol {\tau}}\}&=({\boldsymbol {\tau}}^+ +{\boldsymbol {\tau}}^-), \quad
\lbrack\!\lbrack\,{\boldsymbol{\tau}}\,\rbrack\!\rbrack
={\boldsymbol {\tau}}^+\cdot{\bf n}^+ + {\boldsymbol {\tau}}^-\cdot{\bf n}^- &\mbox{on } E\in
\mathcal{E}^{o}_h, \nonumber
\end{array}
\end{equation}
\begin{equation}\label{ayusob_plenary_av-jump-boundary}
{\lbrack\!\lbrack\, \zeta \,\rbrack\!\rbrack}
=\zeta{\bf n}, \qquad \qquad\qquad
\{{\boldsymbol {\tau}}\} = {\boldsymbol {\tau}} \qquad \mbox{on } E\in \mathcal{E}^{\partial}_h.
\end{equation}
We will also use the notation
$$
(u,w)_{\mathcal{T}_h}=\displaystyle\sum_{T \in \mathcal{T}_h} \int_{T} uw dx \qquad \langle u, w\rangle_{\mathcal{E}_h}=\displaystyle\sum_{E\in \mathcal{E}_h} \int_{E} u w \quad \forall\, u,w,\in V_{h}^{DG}. $$
The approximation to the solution of (\ref{ayusob_plenary_mod0}) reads:
\begin{equation}\label{ayusob_plenary_ayuso_contrib:0}
\mbox{Find }\quad u\in V_{h}^{DG}\quad \mbox{such that }\quad \mathcal{A}_{h}(u,w)=(f,w)_{\mathcal{T}_h}\;, \quad \forall\, w\in V_{h}^{DG}\;,
\end{equation}
with $ \mathcal{A}_{h}(\cdot,\cdot)$  the bilinear form corresponding to  the Interior Penalty (IP) method (see \cite{ayusob_plenary_Arnold82}) defined by:
\begin{equation}\label{ayusob_plenary_ip}
\!\!\!\mathcal{A}_{h}(u,w)\!=\! ( \nabla u,\nabla
w)_{\mathcal{T}_h}\!-\!\langle  {\lbrack\!\lbrack\, u \,\rbrack\!\rbrack}
,\{\nabla w\}\rangle_{\mathcal{E}_h}\! -\langle \{\nabla u\},
{\lbrack\!\lbrack\, w \,\rbrack\!\rbrack}
\rangle_{\mathcal{E}_h}\!+\!\langle S_{h} {\lbrack\!\lbrack\,  u
  \,\rbrack\!\rbrack}
, {\lbrack\!\lbrack\, w \,\rbrack\!\rbrack}
\rangle_{\mathcal{E}_h}\;,
\end{equation}
where $S_{h}=\alpha_{e}\ell_e^{2} h_{e}^{-1}$ with $\alpha_{e}\geq
\alpha^{\ast}>0$ for all $e\in \mathcal{E}_h$,  $h_{e}$ denotes the length of the edge
$e$ in $d=2$ and the diameter of the face $e$ in $d=3$, and $\ell_e=\displaystyle{\max_{T^{+}\cap T^{-}=e}{\{ \ell_{T^{+}}, \ell_{T^{-}}\}}}$, with $\ell_{T^{\pm}}$ being the polynomial degree on $T^{\pm}$.
Following \cite{ayusob_plenary_bcms}, the above IP-biliear form  can be re-written in terms of the weighed residual formulation:
\begin{equation}\label{ayusob_plenary_ayuso_contrib:2}
\mathcal{A}_{h}(u,w)= ( -\Delta u, w)_{\mathcal{T}_h}\!+ \langle
{\lbrack\!\lbrack\, \nabla u \,\rbrack\!\rbrack}
,\{w\}\rangle_{\mathcal{E}^{o}_h} + \langle  {\lbrack\!\lbrack\, u
  \,\rbrack\!\rbrack}
,\,\left(S_h{\lbrack\!\lbrack\, w \,\rbrack\!\rbrack}
 -\{\nabla w\}\right)\rangle_{\mathcal{E}_h}\;. 
\end{equation}
Continuity and Stability can be easily shown in the DG norm or in the
induced $\|\cdot\|_{\mathcal{A}}$-norm, 
provided $\alpha_e\geq \alpha^{\ast}>0$ is taken sufficiently large;
\begin{equation}
\begin{array}{lllll}
&\mbox{ Continuity:   }  \quad\qquad&\mathcal{A}_{h}(u,w) \leq c_c \|u\|_{\mathcal{A}} \|w\|_{\mathcal{A}} \qquad \forall\, u,w\in V_{h}^{DG}\\
&\mbox{ Coercivity:    }\quad \qquad &\mathcal{A}_{h}(u,u) \geq c_s \|u\|_{\mathcal{A}}^{2} \qquad \qquad \quad \forall\, u\in V_{h}^{DG}
\end{array}
\end{equation}

\section{Non-overlapping Domain Decomposition Schwarz
  methods}\label{ayusob_plenary_sec:3}
To define the non-overlapping preconditioners, we need to introduce some further notation.  
We denote by $\mathcal{T}_{S}$ the family of partitions of $\Omega$ into $N$ non-overlapping subdomains $\Omega=\cup_{i=1}^{N} \Omega_i$. Together with $\mathcal{T}_{S}$, we let
 $\mathcal{T}_H$ and $\mathcal{T}_h$ be two families of coarse and fine partitions, respectively, with mesh sizes $H$ and $h$. The three families of
partitions are assumed to be shape-regular and nested:$\mathcal{T}_{S} \subseteq\mathcal{T}_{H} \subseteq\mathcal{T}_{h}$.\\
Similarly as we did for $\mathcal{T}_{h}$ in Section \ref{ayusob_plenary_sec:2}, we define the skeleton and the corresponding sets of internal and boundary edges relative to the subdomain partition.
In particular, for each subdomain $\Omega_i\in\mathcal{T}_{S}$ we define the sets of internal $\mathcal{E}_i^{o}=\{e \in \mathcal{E}_h \,:\, e \subset \Omega_i\}$ and boundary edges $\mathcal{E}_i^{\partial}=\{e \in \mathcal{E}_h \,:\, e \subset \partial \Omega_i \}$, and we set   $\mathcal{E}_i=\mathcal{E}_i^{o}\cup  \mathcal{E}_i^{\partial}$.
Finally, we denote by $\Gamma$ the collection of all interior edges that belong to the skeleton of the subdomain partition;
$$
\varGamma=\bigcup_{i=1}^{N}\varGamma_{i}\;, \quad \mbox{ with  } \quad \varGamma_{i} =\{ e\in\mathcal{E}_h^{o} \,: \, e\subset \partial \Omega_i\}.
$$
The subdomain partition $\mathcal{T}_{S}$ induces a natural space splitting of the $V^{DG}$ finite element space. More precisely,  we have a local finite element subspace associated to each $\Omega_i$  for each 
$i=1,\ldots,S$,  defined by
\begin{equation}
V_h^i=\{ w \in V^{DG}  \,:\,w\equiv 0  \quad \mbox{in } \subset \Omega\smallsetminus \overline{ \Omega_i} \}.
\end{equation}
Let  $\mathcal{I}_{i}^{T}:V_{h}^{i}\longrightarrow V_{h}^{DG}$ be the {\it prolongation} operator, defined as the standard inclusion operator that maps functions of $V_{h}^{i}$ into $V^{DG}_{h}$. We denote by $\mathcal{I}_{i}$ the corresponding {\it restriction} operators  defined (for each $i$) as the transpose of $\mathcal{I}_{i}^{T}$ with respect to the $L^{2}$--inner  product. 
For vector-valued functions $\mathcal{I}_{i}^{T}$ and $\mathcal{I}_{i}$ are defined componentwise. Then the following splitting holds (orthogonal with respect to $L^{2}$-inner product):
\begin{equation}\label{ayusob_plenary_decom:0}
 V_{h}^{DG} =\mathcal{I}_1^T V_h^1\oplus \mathcal{I}_2^T V_h^2\oplus \ldots \oplus \mathcal{I}_{N}^T V_h^{N}\;. 
\end{equation}
\noindent {\sc Local Solvers:} Two types of local solvers have been considered:
\begin{itemize}
\item[{\bf (a).}  ]   {\it $\,\,\,$ Exact local solvers:} Following  \cite{ayusob_plenary_kara0}, the local solvers are defined as the restriction of the discrete bilinear form to the subspace $V_i$. 
\begin{equation}\label{ayusob_plenary_local:kara}
a_i(u_i,w_i)= \mathcal{A}_h (\mathcal{I}_{i}^{T}u_i,\mathcal{I}_{i}^{T}w_i) \qquad \forall\, u_i,w_i \in V^{i}_h
\end{equation}

\item[{\bf (b).}  ]   {\it $\,\,\,$ Inexact local solvers:} Following \cite{ayusob_plenary_paola-blanca1,ayusob_plenary_paola-blanca2} the local solvers are defined as the IP approximation to the original problem (\ref{ayusob_plenary_mod0}) but restricted to the subdomain $\Omega_i$; i.e.,
 \begin{equation}
-\varDelta u^{\ast}_i = f|_{\Omega_i} \quad\mbox{in    }\Omega_i, \qquad u^{\ast}_{i}=0  \quad \mbox{on    }   \partial  \varOmega_i\;.
\end{equation}
Then, the bilinear form can be written as:
\begin{equation}\label{ayusob_plenary_local:inexact}
\widehat{a}_i(u_i,w_i)= ( -\Delta u_i, w_i)_{\mathcal{T}_h\cap
  \Omega_i}\!+ \langle  {\lbrack\!\lbrack\, \nabla u_i \,\rbrack\!\rbrack}
,\{w_i\}\rangle_{\mathcal{E}_i^{o}} + \langle  {\lbrack\!\lbrack\, u_i
  \,\rbrack\!\rbrack}
,S_h{\lbrack\!\lbrack\, w_i \,\rbrack\!\rbrack}
 -\{\nabla w_i\}\rangle_{\mathcal{E}_i}\;,
\end{equation}
where in the above definition, edges on $\mathcal{E}_i^{\partial}$ are
regarded as boundary edges 
(even those $e\in \mathcal{E}_i^{\partial}\smallsetminus \partial \Omega_i$ so that  $e\in \mathcal{E}_h^{o}$) and therefore the trace operators on such edges are defined as in (\ref{ayusob_plenary_av-jump-boundary}).
\end{itemize}
Observe that, in a conforming framework, the definitions given in {\bf
  (a)} and {\bf (b)} would have given rise to exactly the same local
solvers. The difference in the DG context, originates from the
distinct definition of the trace operators on boundary and internal
edges and the fact that $e\in
\mathcal{E}_i^{\partial}\smallsetminus \partial \Omega_i$ is an
interior edge for the global IP method (and so for
(\ref{ayusob_plenary_local:kara})), but a boundary edge for
(\ref{ayusob_plenary_local:inexact}). See
\cite{ayusob_plenary_paola-blanca1,ayusob_plenary_paola-blanca2} for
further details.

Let now $\mathbb{A}$ be the matrix representation of the operator
associated to the global IP method (\ref{ayusob_plenary_ip}), in some
chosen basis (say nodal lagrange basis functions to fix ideas).  We
denote by $\mathbb{A}_i$ and $\widehat{\mathbb{A}}_i$ the matrix
representation (stiffness matrix) of the operators associated to
(\ref{ayusob_plenary_local:kara}) and
(\ref{ayusob_plenary_local:inexact}), respectively. At the algebraic
level, a one-level Additive Schwarz preconditioner is then defined by
$B_{add}^{one}= \sum_{i=1}^{S} \mathbb{I}_{i}^{T} \mathbb{S}_i^{-1}
\mathbb{I}_{i}$ where $\mathbb{I}_{i}$ is the matrix representation of
the restriction operator and $\mathbb{S}_i$ denotes here the matrix
representation of the local solver; and can be chosen to be either
$\mathbb{A}_i$ or $\widehat{\mathbb{A}}_i$.  Notice however, that only
for the choice $\mathbb{S}_i=\mathbb{A}_i$, the resulting one level
additive Schwarz method $B_{add}^{one}$ corresponds to the standard
block jabobi preconditioner for the global stiffness matrix
$\mathbb{A}$. This can be easily checked by noting that the definition
(\ref{ayusob_plenary_local:kara}) gives at the algebraic level
$\mathbb{A}_i= \mathbb{I}_{i}\mathbb{A} \mathbb{I}_{i}^{T}$; that is,
the matrices $\mathbb{A}_i$ are the principal submatrices of
$\mathbb{A}$.
In contrast,   the one level additive Schwarz based on the choice $\mathbb{S}_i=\widehat{\mathbb{A}}_i$ cannot be obtained by starting directly from the algebraic structure of the global matrix $\mathbb{A}$; it would require further modifications of the prolongation and restriction operators.\\

On the other hand, in view of the possibility of considering (at least) these two definitions for the local solvers, a natural question arises. Namely, 
if the inexact local solvers (\ref{ayusob_plenary_local:inexact}) are approximating the original PDE restricted to the subdomain, {\it which continuous problem is approximated by the exact local solvers  (\ref{ayusob_plenary_local:kara}), if any.} By rewriting the bilinear form in the weighted residual formulation one easily obtains:
\begin{equation}\label{ayusob_plenary_local:2}
\begin{array}{llll} 
a_i(u_i,w_i)&= ( -\Delta u_i, w_i)_{\mathcal{T}_h\cap
  \Omega_i}\!+ \langle  {\lbrack\!\lbrack\, \nabla u_i \,\rbrack\!\rbrack}
,\{w_i\}\rangle_{\mathcal{E}_i^{o}} &\\
&\qquad  + \langle  {\lbrack\!\lbrack\, u_i \,\rbrack\!\rbrack}
, \left( S_h{\lbrack\!\lbrack\, w_i \,\rbrack\!\rbrack}
 -\{\nabla w_i\}\right) \rangle_{\mathcal{E}_i^{o}\cup (\mathcal{E}_i^{\partial}\cap \partial\Omega)} &\\
&\qquad +   \langle  \frac{1}{2} \nabla u_i\cdot {\bf n} +S_h u_i  ,   w_i\rangle_{\Gamma_i} - \langle u_i,  \frac{1}{2} \nabla w_i\cdot {\bf n} \rangle_{\Gamma_i} \qquad&
\end{array}
\end{equation}
The terms on the first and second lines are easy to recognize, the first imposes the PDE on each element; the second is the consistency term and the terms in the second line ensure stability and symmetry.  
As  regards those in the last line, the first term is imposing the boundary condition on $\Gamma_i$ (the  part of $\partial\Omega_i\smallsetminus \partial\Omega$). The second term, could be regarded as an artifact to ensure the symmetry of the method. Then, one can write the continuous problem 
\begin{equation}\label{ayusob_plenary_local:cont}
\left\{ \begin{array}{ccccccc} 
-\Delta u^{\ast}_i &=& f|_{\Omega_i} \quad &\mbox{in    }  \Omega_i, \\
 u^{\ast}_{i}&=& 0  \quad & \mbox{on    }   \partial \varOmega_i\cap  \partial \varOmega\;,\\
 \frac{1}{2}\frac{\partial u^{\ast}_i}{\partial n_i} +S_h u^{\ast}_i &=&0 \quad & \mbox{on    }  \Gamma_i\;.
\end{array}\right.
\end{equation}
This implies that the exact local solvers for the IP method (and in
general for most DG methods) are approximating the original problem
but with transmission Robin conditions. And as $h\to 0$ the method
enforces $u^{\ast}_i=0$ on $\Gamma_i$.  
Whether such interface
boundary conditions are optimal or could be further tuned to improve the
convergence properties of the classical Schwarz methods is a subject
of current research. Optimization of the Schwarz methods with respect
to the interface boundary conditions has been recently studied in 
\cite{ayusob_plenary_soheil}. The final ingredient needed to define the two-level Schwarz method is the coarse solver.

\noindent {\sc Coarse solver:} Let $V_c:=V_{H}^{DG}$ be the coarse space and let  $a_{c}:V_{c}\times V_{c} \longrightarrow \mathbb{R}$ be the coarse solver defined by \cite{ayusob_plenary_kara0,ayusob_plenary_paola-blanca1,ayusob_plenary_paola-blanca2}:
\begin{equation}\label{ayusob_plenary_coarse}
a_c(u_c,w_c)=\mathcal{A}_{h} (\mathcal{I}_{c}^{T}u_c,\mathcal{I}_{c}^{T}w_c) \qquad \forall\, u_c,w_c \in V_{c}
\end{equation}
where $\mathcal{I}_{c}^{T}:V_c\longrightarrow V_h^{DG}$ is the prolongation operator, defined as the standard inclusion. Notice that with this definition, the corresponding matrices do indeed satisfy the Galerkin property: $\mathbb{A} =\mathbb{I}_{c}^{T} \mathbb{A}_c  \mathbb{I}_{c}$, but should be noted that  unlike in a conforming framework $a_c(u_c,w_c)\ne  \mathcal{A}_{H}(u_c,w_c)$. A two level Schwarz preconditioner can then be defined:
\begin{equation}\label{ayusob_plenary_2level}
\mathbb{B}_{add}= \sum_{i=1}^{S}  \mathbb{I}_{i}^{T} \mathbb{S}_i^{-1}  \mathbb{I}_{i} + \mathbb{I}_{c}^{T} \mathbb{A}^{-1}_c  \mathbb{I}_{c}
\end{equation}
It is also possible to define the coarse solver as IP approximation (with the partition $\mathcal{T}_{H}$ and the coarse space $V_c$) to the orginal problem (i.e., as $\mathcal{A}_{H}(u_c,w_c)$). However with such definition, the Galerkin property is lost and in order to ensure scalability of the resulting two level Schwarz preconditioner, more sophisticated prologation and  restriction operators are required \cite{ayusob_plenary_bren10}. \\

Let now $B^{-1}$ denote the inverse operator associated to the two level preconditioner (\ref{ayusob_plenary_2level}). To analyze the convergence properties of the resulting preconditioner one needs to characterize the dependence of the constants $C_1$ and $C_0$ in
\begin{equation}\label{ayusob_plenary_espectral:0}
C_1 \mathcal{A}_{h}(w,w) \leq (B^{-1}w,w) \leq C_0^{2} \mathcal{A}_{h}(w,w)\qquad \forall\, w\in V_{h}^{DG}
 \end{equation}
 The condition number of the preconditioned matrix
 $\mathbb{B}\mathbb{A}$ is then 
$C_0^{2}/C_1$. The proof of (\ref{ayusob_plenary_espectral:0}) is often guided by
 Lions lemma (for a proof see \cite{ayusob_plenary_xu92}, 
\cite{ayusob_plenary_Widlund92},
\cite[Lemma 2.4]{ayusob_plenary_XZ}), which tells that the preconditioner can be written as
\begin{equation}
(B^{-1}w,w):= \inf_{\begin{array}{cc}
&w_i \in V^{i} \\
&w_c+\sum_{i} w_i=w
\end{array}} \left( a_c( w_c,w_c)+ \sum_i \mathcal{R}_i(w_i, w_i) \right),
\end{equation}
where we have denoted by $ \mathcal{R}_i(\cdot,\cdot)$ the {\it
  approximate (or exact) subspace solver} on $V^{i}$. 
\section{Ficticious Space and Auxiliary Space Methods}\label{ayusob_plenary_sec:4}
Ficticious Space Lemma was originally introduced by Nepomnyaschikh in
\cite{ayusob_plenary_NEP1991}, and further used for developing and analyzing
multilevel preconditioners for nonconforming approximations in
\cite{ayusob_plenary_Oswald96}  and for conforming methods with nonconforming meshes in \cite{ayusob_plenary_JXU96}.  There
are two main ingredients to construct a fictitious space
preconditioner for the operator $A: V_{h}^{DG}\longrightarrow V_{h}^{DG}$
associated to the bilinear form (\ref{ayusob_plenary_ip}).
\begin{enumerate}
\item[(1)] A fictitious space $\overline{V}$, and an symmetric positive
  definite operator 
  $\overline{A}:\overline{V}\longrightarrow \overline{V}$  associated with some  $\overline{\mathcal{A}}(\cdot,\cdot):\overline{V}\times \overline{V}\longrightarrow \mathbb{R}$.
  \item[(2)] A continuous, linear and surjective mapping $\Pi: \overline{V}\to V_{h}^{DG}$
\end{enumerate}
The fictitious  space preconditioner $B$ is then defined as
\begin{equation}\label{ayusob_plenary_eqn:Bdef}
  B=\Pi\circ{\overline{A}}^{-1}\circ \Pi^*: V_{h}^{DG} \to V_{h}^{DG}.
\end{equation}
The convergence properties of  the preconditioner
 $B$ depend on the choice of the fictitious space $\overline{V}$ and ficticious operator
$\overline{A}$. Typically, one chooses a fictitious pair $(\overline{V},\overline{A})$ for which it is simpler to construct a preconditioner.
 The analysis of such
methods is done via the \emph{Fictitious space lemma}~\cite{ayusob_plenary_NEP1991},
which states that if $\Pi$ has a bounded (in energy norm) right inverse
and is stable in $\overline{A}$ norm, then $B$ is equivalent to $A$ (in the sense that they satisfy a corresponding (\ref{ayusob_plenary_espectral:0})) with
constants of equivalence  ($C_1$ and $C_0^{2}$) depending on the stability and invertibility
of $\Pi$. 
The auxiliary space idea, comes from the observation (see~\cite{ayusob_plenary_JXU96}) that  a \emph{surjective} $\Pi$
is easy to construct for the choice $\overline{V}=V_{h}^{DG}\times W$ for some space $W$ (the factor $V_{h}^{DG}$ in
the product plays a crucial role). 

One natural approach in constructing such preconditioners for DG
discretizations is via subspace splitting which uses the corresponding
conforming space as the component $W$; that is
$\overline{V}=V_{h}^{DG}\times V_{\tilde{h}}^{\rm{conf}}$, with
$W:=V_{\tilde{h}}^{\rm{conf}}$ denoting the conforming finite element
space with $ \tilde{h}$ chosen $\tilde{h}\geq h$. This is natural
because one expects that the smooth error (with small energy) is in
this space. Then, for the auxiliary preconditioner $\overline{A}^{-1}$
one can choose his favourite solver in
$V_{\tilde{h}}^{\rm{conf}}$. Preconditioners based on such splittings
are found in~\cite{ayusob_plenary_dobrev} and
~\cite{ayusob_plenary_dahmen1}, and more recently in
\cite{ayusob_plenary_kolja1,ayusob_plenary_luca}.  Two-level methods
based on three different splittings of the DG space are given
in~\cite{ayusob_plenary_dobrev}. In~\cite{ayusob_plenary_dahmen1}, an
auxiliary space preconditioner is proposed (and analyzed) for IP
discretizations with non-conforming meshes and hanging nodes. This
auxiliary space approach has been recently extended and used for
designing multilevel preconditioners in \cite{ayusob_plenary_kolja1}
for the IP method with arbitrary polynomial degree. The results
from~\cite{ayusob_plenary_kolja1} are further used for constructing a
BDDC preconditioner for such discretizations in
\cite{ayusob_plenary_luca}.

We wish to point out that for the IP method such decompositions were
already known in the area of adaptivity and a posteriori error
analysis for DG methods.  The following important decomposition is
implicitly contained in In the seminal
work~\cite{ayusob_plenary_kara-apost1}:
\begin{equation} 
V_h^{DG} = V_h^{\rm{conf}} \oplus E_h,
\end{equation} 
where $E_h=(V_h^{\rm{conf}})^{\perp}$ refers to the complementary space of $V_h^{\rm{conf}}$
 in $V^{DG}$ (orthogonal with respect to the corresponding energy inner product).
In fact,  an explicit construction of an
interpolation operator $I_h: V_h^{DG}\longrightarrow V_h^{\rm{conf}}$ is provided,
on simplicial meshes, even in case of hanging nodes, which is stable in the energy norm,
and therefore can be used as a component in constructing a stable
surjective $\Pi$ in the design of an auxiliary space preconditioner. 

The analysis of the auxiliary space preconditioners using the conforming method as a
component of the space decomposition is carried out in a standard fashion
by introducing stable and accurate interpolation operators (see
e.g. ~\cite{ayusob_plenary_dahmen1} or \cite{ayusob_plenary_dobrev} for such
constructions). Alternatively, at least for the $h$-version, one may adapt and use the framework developed
in~\cite{ayusob_plenary_kara-apost1} to analyse the properties of these
preconditioners.
\section{Orthogonal space splittings in a nutshell}\label{ayusob_plenary_sec:5}
The approach we present now has been developed in \cite{ayusob_plenary_az} for developing uniform solvers for the family of IP discretizations, including non-symmetric schemes. 
It could be seen as a clever change of basis which allows
for special decompositions of the DG space. The ideas work in dimensions $d=2,3$ and are based on a natural splitting of the linear DG  FE space on simplicial meshes with no-hanging nodes. Therefore, in all what follows $V^{DG}$ stands for the linear approximation space; i.e., $\ell=1$.
Furthermore, to ease the presentation,  we drop the subindex $h$ from the finite element space and the bilinear form, so $\mathcal{A}(\cdot,\cdot)= \mathcal{A}_{h}(\cdot,\cdot)$. For multilevel considerations see for instance \cite{ayusob_plenary_jump}.
To introduce the space splitting we  first introduce some notation. 

\begin{figure}[!htp]
        \begin{center}
        \vskip -0.1cm
                          \includegraphics[width=4.1cm]{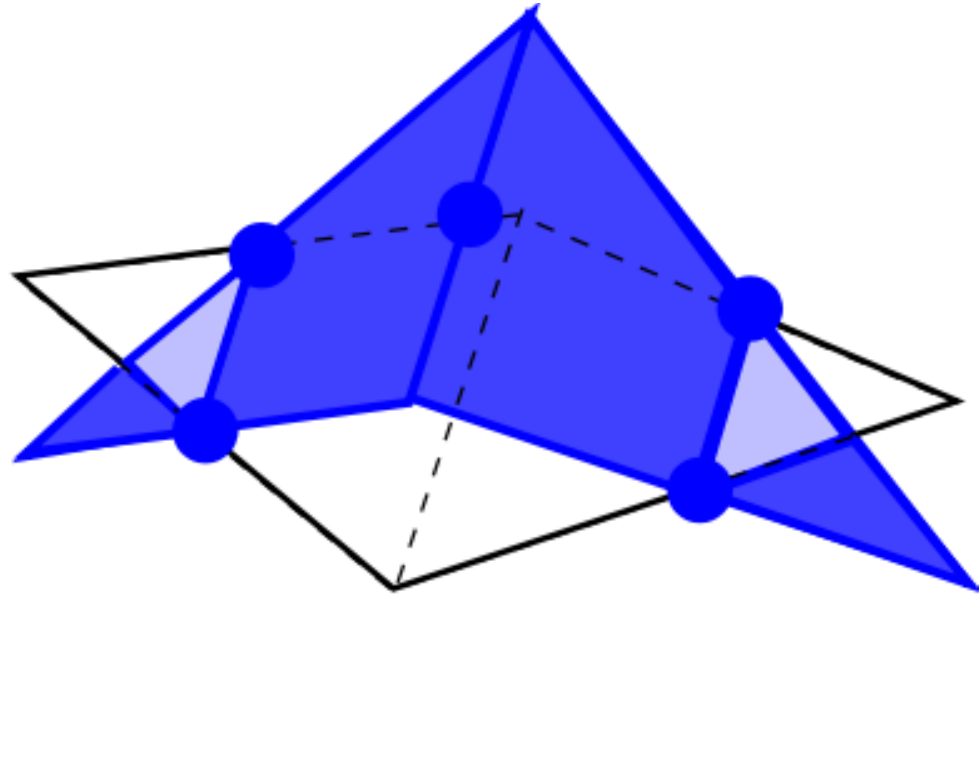}
        \hskip 1.2cm
          \includegraphics[width=4.1cm]{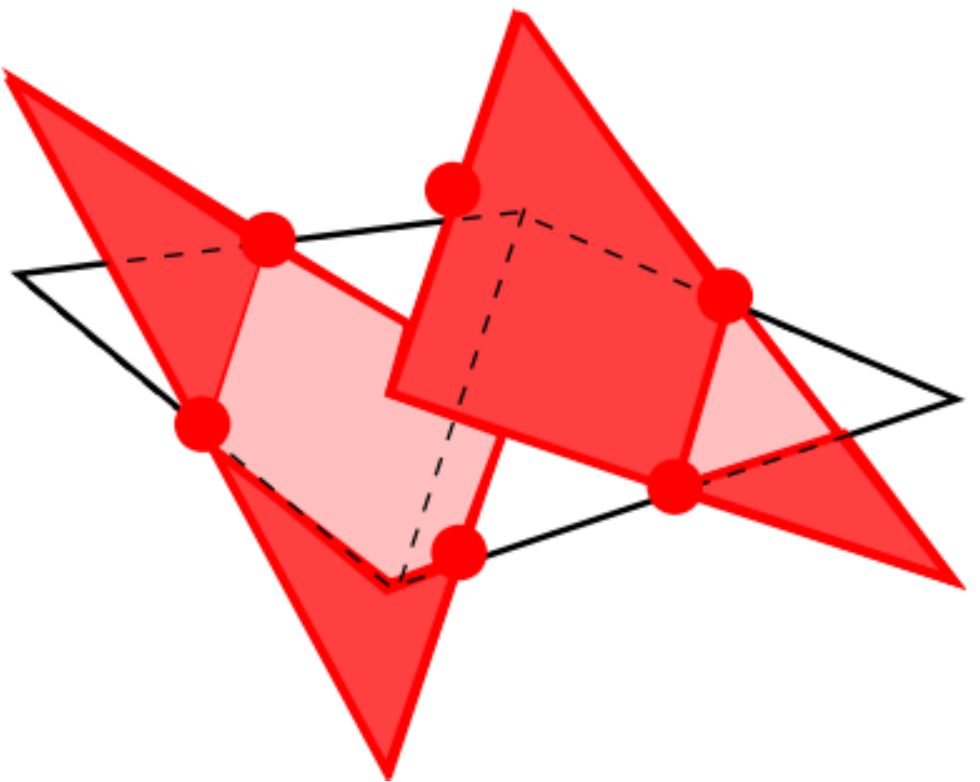}    
                    \end{center}
                    \vskip -0.1cm
             \caption{Basis functions (associated to an edge) for the Crouziex Raviart space (left figure) and the $\mathcal{Z}$ space (right figure)}
        \label{ayusob_plenary_fig0}
       \end{figure}

Together with the IP  bilinear form $\mathcal{A}(\cdot,\cdot)$, we also consider  the bilinear form that results by computing all the integrals in (\ref{ayusob_plenary_ip}) with the mid-point quadrature rule, known as weakly penalized or IP-0 method:
\begin{equation}\label{ayusob_plenary_ip0}
\mathcal{A}_0(u,w)= (-\Delta u, w)_{\mathcal{T}_h} +\langle
{\lbrack\!\lbrack\, \nabla u \,\rbrack\!\rbrack}
, \{w\}\rangle_{\mathcal{E}^{o}_h} +\langle
\mathcal{P}^{0}_{E}({\lbrack\!\lbrack\,  u \,\rbrack\!\rbrack}
), S_h {\lbrack\!\lbrack\,  w \,\rbrack\!\rbrack}
- \{\nabla w\}\rangle_{\mathcal{E}_h}\;,
\end{equation}
where, for each $e\in \mathcal{E}_h$, let $\mathcal{P}_{e}^{0}:L^{2}(e)\longrightarrow \mathbb{P}^{0}(e)$ is  the $L^{2}$-orthogonal projection onto the
constants on that edge defined by:
\begin{equation}\label{ayusob_plenary_proj0}
\mathcal{P}_{e}^{0}(u):= \frac{1}{|e|} \int_{e} u,  \qquad \forall\, u\in L^{2}(e).
\end{equation}
 We define the following two subspaces of $V^{DG}$
 \begin{eqnarray}
 V^{CR}&:=\{ v\in V^{DG} \,\,\, :\quad
 \mathcal{P}_{e}^{0}({\lbrack\!\lbrack\, v \,\rbrack\!\rbrack}
)=0 \,\,\, \forall\, e\in \mathcal{E}_h^{o}\} \label{ayusob_plenary_defCR}\\
  \mathcal{Z}&:=\{ z\in V^{DG} \,\,\, :\quad \mathcal{P}_{e}^{0}(\{z\})=0 \,\,\, \forall\, e\in \mathcal{E}_h\} \label{ayusob_plenary_defZ}
\end{eqnarray}
The first one is the well known lowest order Crouziex-Raviart finite element space. 
The above subspaces can be seen to be complementary to each other, and in fact it is easy to prove that
\begin{equation}\label{ayusob_plenary_decom}
V^{DG}=V^{CR} \oplus \mathcal{Z}\;.
\end{equation}
Notice that the explicit characterization of the subspaces allows to provide basis for both spaces. (See Fig. \ref{ayusob_plenary_fig0}).
      
A key property satisfied by the space decomposition (\ref{ayusob_plenary_decom}) is that the two subpaces are orthogonal in the enegy norm defined by $\mathcal{A}_0(\cdot,\cdot)$. In fact it can be easily shown using (\ref{ayusob_plenary_ip0}) and the definition of the spaces (\ref{ayusob_plenary_defCR}) and (\ref{ayusob_plenary_defZ})  that
\begin{equation}\label{ayusob_plenary_orto}
\mathcal{A}_0(v,z)=\mathcal{A}_0(z,v)=0 \qquad \forall\, v\in V^{CR}\;, \,\, z\in \mathcal{Z}\;.
\end{equation}
This already suggest that by perfoming a {\it change of basis} of the standard Lagrange basis for $V^{DG}$ to the ones in $V^{CR}$ and $\mathcal{Z}$, the stiffness matrix representation of $\mathbb{A}_0$  in the new basis have a block diagonal structure.  Therefore, for the IP-0 method the following algorithm is an exact solver:\\

\noindent {\it Algorithm 1:} Let $u_0$ be a given initial guess. For $k\ge 0$, and given
$u_k=z_k+v_k$, the next iterate
$u_{k+1}=z_{k+1}+v_{k+1}$ is defined via the two steps:
\begin{enumerate}
\item[1.] Solve $\mathcal{A}_0(z_{k+1},\psi^{z}) = (f,\psi^{z})_{\mathcal{T}_h}
 \quad \forall\, \psi^{z}\in \mathcal{Z}$.
\item[2.] Solve $\mathcal{A}_0(v_{k+1},\varphi) = (f,\varphi)_{\mathcal{T}_h}
\quad \forall\, \varphi\in V^{CR}$.
\end{enumerate}
Notice that algorithm 1 requires two solutions of smaller
 problems: one solution in $\mathcal{Z}$-space (step $1$ of the
algorithm~1), and one solution in $V^{CR}$-space (step $2$ of
algorithm~1).  
As we show next,  the solution of 
 the subproblems on $\mathcal{Z}$ and on $V^{CR}$ can be done efficiently.\\

\noindent {\sc Solution in the $\mathcal{Z}$-space:} The functions in $\mathcal{Z}$ have non-zero jump on every edge, which suggest the high oscillatory nature of its functions.  Using the definition of the space, the following useful property (Poincare-type inequality) can be shown:
\begin{lemma}\label{ayusob_plenary_le:poincZ}
Let $\mathcal{Z}$ be the space defined in (\ref{ayusob_plenary_defZ}). 
$$
h^{-2}\| z\|_{0,\mathcal{T}_h}^{2}\lesssim \mathcal{A}_0(z,z)\lesssim
h^{-2}\| z\|_{0,\mathcal{T}_h}^{2}\;, \qquad \forall\, z\in
\mathcal{Z}
$$
\end{lemma}
By virtue of this lemma it follows that the condition number (denoted
by $\kappa$) of the block matrix associated to the restriction of
$\mathcal{A}_0(\cdot,\cdot)$ to the subspace $\mathcal{Z}$, say
$\mathbb{A}^{zz}_0$, satisfies $\kappa (\mathbb{A}_0^{zz})=O(1)$ and it
is independent of the mesh size. Therefore, efficient solver for the
problem in $\mathcal{Z}$ is the Conjugate
Gradient (CG) method with a simple diagonal preconditioner.\\

\noindent  {\sc Solution in $V^{CR}$:} The restriction of
$\mathcal{A}_0(\cdot,\cdot)$ to the $V^{CR}$ sub-space gives the
well-known 
Crouziex-Raviart approximation method for (\ref{ayusob_plenary_mod0}) ;
 \begin{equation}\label{ayusob_plenary_a00}
 \mathcal{A}_0(v,\varphi)=(\nabla v,\nabla \varphi)_{\mathcal{T}_h}=\sum_{T\in \mathcal{T}_h}(\nabla v,\nabla\varphi)_{T} \quad \forall\, \,v\, ,\varphi \in V^{CR}\;,
\end{equation}
Therefore, it is enough to resort to any of the solvers that have been already developed, for instance \cite{ayusob_plenary_brennerNC,ayusob_plenary_Oswald96,ayusob_plenary_sarkisNC0}.\\

So far, an exact solver has been constructed in a simple and clean way for the IP-0 method.
A last ingredient is needed to provide uniformly convergent solvers for the IP method (\ref{ayusob_plenary_ip}) and it is formulated in next Lemma:
\begin{lemma}\label{ayusob_plenary_ayuso_contrib:4}
Let $\mathcal{A}(\cdot,\cdot)$  and  $\mathcal{A}_{0}(\cdot,\cdot)$ be
the bilinear forms of the IIPG method defined in (\ref{ayusob_plenary_ip}) and (\ref{ayusob_plenary_ip0}). Then, there exist $c_2>0$ depending only on the shape regularity of $\mathcal{T}_h$ and $c_0>0$ depending also on the penalty parameter $\alpha$ such that 
    \end{lemma}
\begin{equation}\label{ayusob_plenary_equivA:A0}
 c_2 \mathcal{A}_{0} (u,u)\lesssim \mathcal{A}(u,u)\leq c_0 \mathcal{A}_{0} (u,u)\quad \forall u\in V^{DG}.
\end{equation}
The above result establishes the \emph{spectral equivalence} between $\mathcal{A}_0(\cdot ,\cdot)$ and $\mathcal{A}(\cdot,\cdot)$. Therefore, in terms of solution techniques, a uniform preconditioner for the IP-0 method, already provides a uniform preconditioner for the IP method.\\

These ideas and new framework,  have  been already extended and adapted for designing and analyzing solvers for other problems:\\

\indent $\bullet$ In  \cite{ayusob_plenary_jump} the case of second order elliptic problems with large {\it jumps in the diffusion coefficient} is considered. In a first step, the space splitting (\ref{ayusob_plenary_decom}) needs to be modified to account for the jumps in the coefficient, while still being orthogonal with respect to the corresponding $\mathcal{A}_0(\cdot,\cdot)$-induced norm.  The choice of a robust method for approximating the continuous problem (definition of the relevant $\mathcal{A}(\cdot,\cdot)$ bilinear form) allows to guarantee that the corresponding spectral equivalence property (\ref{ayusob_plenary_equivA:A0}) holds with constants $c_0,c_2$ independent of the mesh size and the {\it jumping coefficient}.\\
\indent $\bullet$ In \cite{ayusob_plenary_AyusoB_GeorgievI_KrausJ_ZikatanovL-2009aa} efficient
solvers are analyzed for IP approximations of {\it linear elasticity
  problems}, considering all cases: the pure displacement, the mixed
and the traction free problems.
The last two cases pose some extra pitfalls in the analysis  since the spectral equivalence property (\ref{ayusob_plenary_equivA:A0}) does not hold in those cases. In spite of that, the ideas can still be  used to construct block preconditioners (guided by the algebraic structure of $\mathcal{A}_0(\cdot,\cdot)$ due to the orthogonality) and prove uniform convergence.\\
\indent $\bullet$ In \cite{ayusob_plenary_ariel00} it is shown how to construct an efficient
solver for the solution of the linear system that arise from a DG
discretization of a convection-diffusion problem, in the convection
dominated regime. The problem is relevant in semiconductor
applications. In this case, the original method is a non-symmetric
exponentially fitted IP weakly-penalized.
\section*{Acknowledgments}
B. Ayuso de Dios thanks R. Hiptmair (ETH) for raising up the issue on the use of Auxiliary space techniques for preconditioning DG methods, at the DD21 meeting.
First author has been partially supported by MINECO grant MTM2011-27739-C04-04 and GENCAT 2009SGR-345.


\bibliographystyle{spmpsci}
\bibliography{arxiv_ayusob_plenary} 


\end{document}